\documentclass{amsart}
\usepackage{amsfonts,latexsym}
\newtheorem{theorem}{Theorem}

\newtheorem{proposition}{Proposition}
\newtheorem{corollary}{Corollary}

\usepackage{amsfonts}

\newcommand{\beq}{\begin{equation}}
\newcommand{\eeq}{\end{equation}}

\begin{document}

\title[Coxeter group actions on the complement of hyperplanes]{Coxeter group actions on the complement of hyperplanes and special involutions}
       \author{G. Felder}
       \address{Department of Mathematics,
          ETH, Zurich,
Switzerland}
       \email{felder@math.ethz.ch}

       \author{A.P. Veselov}
       \address{Department of Mathematical Sciences,
        Loughborough University, Loughborough,
        Leicestershire, LE11 3TU, UK and
       Landau Institute for Theoretical Physics, Moscow, Russia}
       \email{A.P.Veselov@lboro.ac.uk}

\maketitle

%\bigskip

\bigskip

{\small  {\bf Abstract.} We consider both standard and twisted
action of a (real) Coxeter group $G$ on the complement $\mathcal M_G$
to the complexified reflection hyperplanes by combining the
reflections with complex conjugation. We introduce a natural geometric class of special involutions in $G$ and give explicit formulae
which describe both actions on the total cohomology $H^*(\mathcal M_G,
{\mathbb C})$ in terms of these involutions. 
As a corollary we prove that the corresponding twisted
representation is regular only for the symmetric group $S_n$, the Weyl groups of type $D_{2m+1}$, $E_6$ and dihedral groups $I_2 (2k+1).$ We discuss also the relations with the cohomology of Brieskorn's braid groups.}

\bigskip

\section*{Introduction.}

In 1969 V.I. Arnol'd \cite{Arnold} computed the cohomology algebra
of the configuration space $\mathcal M_n$ of $n$ distinct points of
the complex plane. This remarkable short paper was the starting
point of the active research in this area of mathematics on the
crossroad of algebra, geometry and combinatorics.

In particular, Brieskorn \cite{Brieskorn} generalised Arnol'd's
results to arbitrary irreducible Coxeter group $G$ and showed that the
Poincar\'e polynomial of the complement $\mathcal M_G$ to the
complexified reflection hyperplanes has the form
$$P(\mathcal M_G,t) = (1 + m_1 t)\cdots(1 + m_n t),$$ where $m_i = d_i -1$ are the
exponents of the Coxeter group $G,$ $d_i$ are the degrees of the generators in the algebra of $G$-invariants. Since the product $(1 +
m_1)\cdots(1 + m_n) = d_1\cdots d_n= |G|$ is known to be the order of the group $G$,
it is tempting to suggest that the total cohomology space
$H^*(\mathcal M_G) = H^*(\mathcal M_G, {\mathbb C})$ is the regular
representation with respect to the natural action of $G$ on
$\mathcal M_G.$ However this turns out not to be true already for the
symmetric group $G = S_n$, as it was shown by Lehrer
\cite{L1}, although not far from being true.

The starting point of this work was the following observation. Let
us consider another action of $G$ on $\mathcal M_G$ using the
anti-holomorphic extensions of the reflections from the real to
the complexified space. In other words we combine the usual action of
the reflections with complex conjugation.  The claim is that, for
the symmetric group $G = S_n$ with this new action, $H^*(\mathcal M_G,
{\mathbb C})$ is indeed a regular representation. We were not able to
find this result in the literature, although the proof can be easily
derived from known results. While looking for the simplest
explanation of this fact, we have found a simple universal way to
investigate this representation for any Coxeter group, both in the
standard and twisted case.

The main result of this paper is the following two explicit formulae which
show how far are the corresponding representations from the
regular and describe them as virtual representations.

Let $G$ be a finite group generated by reflections in a real Euclidean space $V$ of dimension 
$n,$  $\mathcal M_G$ be the complement to the complexified reflection hyperplanes in $V\otimes {\mathbb C}.$ Let us denote the  twisted representation of
$G$ on $H^*(\mathcal M_G, {\mathbb C})$ by $H_{\epsilon}^*(\mathcal M_G)$.
Alternatively we can write
\begin{equation}
\label{E} H_{\epsilon}^*(\mathcal M_G) = \bigoplus_{k = 0}^{n}
\epsilon^k \otimes H^k(\mathcal M_G, {\mathbb C}),
\end{equation}
where $\epsilon$ is the alternating representation of $G.$

We claim that for all Coxeter groups $G$ with the standard action
on $\mathcal M_G$
\begin{equation}
\label{F1} H^*(\mathcal M_G) = \sum_{\sigma \in X_G} (2 Ind_{<\sigma>}^G(1) - \rho)
\end{equation}
and for the twisted action
\begin{equation}
\label{F2} H_{\epsilon}^*(\mathcal M_G) = \sum_{\sigma \in
X^{\epsilon}_G} (2 Ind_{<\sigma>}^G(1) - \rho).
\end{equation}
Here $\rho$ is the regular representation of $G$, $X_G$ is a
special set of involutions in $G$ (more precisely, conjugacy
classes of involutions), $X^{\epsilon}_G$ is the subset of $X_G$
consisting of even involutions, $Ind_{<\sigma>}^G(1)$ is the
induced representation from the trivial representation of the
subgroup generated by $\sigma.$ 

We should mention that for the standard action in most of the
cases our result can be derived from the results of the papers
\cite{FJ1, FJ2, L1, L2}, so our main contribution here is the following
universal geometric description of the set $X_G.$

Namely, consider any involution $\sigma$ in the geometric realisation
of $G$ as a group of orthogonal transformations of a Euclidean
space $V$  and the corresponding
splitting $V = V_1 \oplus V_2$ so that $V_1=V^-(\sigma)$ and $V_2=V^+(\sigma)$ are the
eigenspaces of $\sigma$ with eigenvalues $-1$ and $1$, respectively. 
%%% I exchanged 1 and -1
Let $R_1$
and $R_2$ be the intersections of the root system $R$ of group $G$
with these subspaces, and $G_1$, $G_2$ be the corresponding
Coxeter subgroups of $G.$ We call the involution $\sigma$ {\it
special} if for any root $a\in R$ at least
one of the projections of $a$ onto $V_1$ and $V_2$ is proportional to a
root from $R_1$ or $R_2$.  In particular the identity and any simple reflection
are always special involutions. 

We denote by $X_G$ the set of all conjugacy classes of special involutions
in $G.$ By choosing a representative in each class we can realise $X_G$ as
a special set of involutions in $G.$

For the symmetric group $S_n,$ Weyl groups of type $D_{2m+1}, E_6$ and dihedral groups $I_2(2k+1)$  the set $X_G$ consists of two elements: the identity and 
the class of a simple reflection $s$, so we have in that case
$$H^*(\mathcal M_G) = 2 Ind_{<s>}^G(1)$$
and
$$H_{\epsilon}^*(\mathcal M_G) = \rho.$$

For the Coxeter groups of type $B_n$ (or $C_n$) the set $X_G$
consists of $2n$ involutions which in the geometric realisation have the
form $\sigma_k = P_{12} \oplus (-I_k) \oplus I_{n-k-2}, k =
0,1,\dots, n-2$ and $\tau_l = (-I_l) \oplus I_{n-l}, l = 0, 1,\dots,n.$ 
Here $P_{12}$ denotes the permutation of
the first two coordinates and $I_k$ is the $k \times k$ identity matrix.

In the case of Weyl groups of type $D_{2m}$, $E_7$, $E_8$ and for the 
icosahedral groups $H_3$ and $H_4$, $X_G$ consists of 4 involutions: 
$\pm Id$ and $\pm s,$ $s$ is a simple reflection.

For $F_4$ we have 8 involutions which in the standard geometric
realisation \cite{Bourbaki} have the form $\pm s$ with $s$
representing two different conjugacy classes of simple
reflections, $diag(-1,-1,1,1),  P_{12} \oplus diag(-1,1)$ and the
centre $\pm Id$.

Finally, for the dihedral groups $I_2(2k)$  $X_G$ consists of 4 elements:
two conjugacy classes of simple reflections and the central
elements $\pm Id$.

%new
In the Appendix we list the graphs of all equivalence classes of nontrivial special involutions using Richardson's description \cite{R}.

As a corollary we have that the twisted action of an irreducible Coxeter group $G$ on 
$H^*(\mathcal M_G)$ is a regular representation only for the symmetric group $S_n$,
the Weyl groups of type $D_{2m+1}$, $E_6$
and dihedral groups $I_2 (2k+1)$. Note that besides the one-dimensional Coxeter group of type $A_1$
this is the list of Coxeter groups with trivial centre.

%%%
Our result can be reformulated in terms of the decomposition of the
cohomology into irreducible $G$-modules: the multiplicity $m(W)$
of the irreducible $G$-module $W$ in the decomposition of $H^*(\mathcal M_G)$ is
\begin{equation}
\label{mW}
m(W) = \sum_{\sigma\in X_G}(\mathrm{dim}\,W^+(\sigma)-\mathrm{dim}\,W^-(\sigma)),
\end{equation}
where $W^\pm(\sigma)=\{v\in W\colon \sigma v=\pm v\}$ 
are the eigenspaces of the involution $\sigma.$ 
For the twisted action, $X_G$ is replaced by $X^\epsilon_G$.

In particular, the trivial representation has the multiplicity
%\begin{equation}
%\label{m1}
$$
m(1) = |X_G|,
$$
%\end{equation}
where $|X|$ denotes the number of elements in the set $X$. 
This gives us a topological interpretation 
of $|X_G|$ as the total Betti number of the corresponding quotient
space $\Sigma_G = \mathcal M_G / G.$
More precisely we show that the Poincar\'e polynomial $P (\Sigma_G, t)$
has the form
\begin{equation}\label{PoinP}
P (\Sigma_G, t) = \sum_{\sigma\in X_G} t^{\mathrm{dim}\,V^-(\sigma)},
\end{equation}
where $V$ is the geometric representation of $G$ and $V^-(\sigma)$ is the $(-1)$-eigenspace
of $\sigma$ in this representation.
Notice that, according to a classical result due to Brieskorn and Deligne, $\Sigma_G$ is the Eilenberg--Mac Lane space $K(\pi,1)$ for the corresponding generalised braid group $\pi = B_G,$ so the last formula also can be interpreted in terms of the (rational) cohomology of $B_G$ (see \cite {Brieskorn}).

Another interesting corollary of our results is that the multiplicity $m(\epsilon)$ of the alternating representation in $H^*(\mathcal M_G)$ is zero. Due to (\ref{mW}) this is equivalent to the fact that the numbers of even and odd classes in $X_G$ are equal. For the group $S_n$ this was conjectured (in different terms) by Stanley \cite{Stanley} and proved by Hanlon in \cite{Hanlon}. For the general
Coxeter groups it was first proved by Lehrer \cite{L3}.

Our approach is geometrical and based on the (generalised) Lefschetz
fixed point formula which says that the Lefschetz numbers for the
action of the finite groups are equal to the Euler characteristic
of the corresponding fixed sets (see \cite{B}). We should mention
that the idea to use the Lefschetz fixed point formula is not new 
in this area (see e.g. \cite{CT}), however in this form it seems to be not explored
before.

\section*{Characters and Lefschetz fixed point formula.}

Consider first the standard action of $G$ on $\mathcal M_G$. The
character of an element $g \in G$ in the corresponding
representation $H^*(\mathcal M_G)$ is $$\chi(g) = \sum_{k=0}^n tr
g_k^*,$$ where $g_k^*$ is the action of $g$ on $H^k(\mathcal M_G).$
Now replace the action of $g$ by its composition with the complex
conjugation which we denote as $\bar g.$ Then because complex
conjugation is acting as $(-1)^k$ on the $k$-th cohomology of
$\mathcal M_G$ we have $$\chi(g) = \sum_{k=0}^n (-1)^k tr \bar
g_k^*,$$ which by definition is the {\it Lefschetz number} $L(\bar
g)$ of the map $\bar g.$ Now we can apply the Lefschetz fixed
point formula \cite{B} which says that this number is equal to the
Euler characteristic $\chi(F_{\bar g})$ of the fixed set $F_{\bar
g} = \{ z \in \mathcal M_G: \bar g(z) = z \}.$

In the twisted case all the even elements of $G$ act in the
standard way but the action of all odd elements are twisted and
the corresponding character of an odd $g \in G$ in
$H_{\epsilon}^*(\mathcal M_G)$ is $$\chi_{\epsilon}(g) = \sum_{k=0}^n
(-1)^k tr g_k^* = L(g)$$ is the Lefschetz number of the standard
action of $g$ on $\mathcal M_G.$

\begin{proposition} \label{l}
\begin{enumerate}
\item[(i)]
The character of an element $g$ from $G$ for the standard action
of $G$ on $H^*(\mathcal M_G)$ is equal to the Euler characteristic of
the fixed set $F_{\bar g}.$
\item[(ii)]
The character of an odd element $g$ from $G$ in
$H_{\epsilon}^*(\mathcal M_G)$ for the twisted action of $G$ on
$\mathcal M_G$ is equal to the Euler characteristic of the fixed set
$F_{g}.$
\item [(iii)]The fixed set $F_{g}$ is empty unless $g$ is identity, so the
character of any odd element in the twisted representation
$H_{\epsilon}^*(\mathcal M_G)$ is zero. The fixed set $F_{\bar g}$ is
empty unless $g$ is an involution, so if the order of $g$ is
bigger than 2 then $\chi(g) = 0.$
\end{enumerate}
\end{proposition}

Only the last part needs a proof.  A fixed point of $g$
corresponds to an eigenvector of $g$ in the geometric realisation with
the eigenvalue $1,$ which can be chosen real (recall that $G$ is a
real Coxeter group). But on the reals $G$ acts freely on the
set of Coxeter chambers.
% I changed the preceding sentence
 Thus $F_g$
is empty unless $g$ is identity. Now assume that $\bar g$ has a
fixed point $z \in \mathcal M_G$, then obviously $h = \bar g^2 = g^2$
must also have a fixed point, which means that $g^2$ must be
identity.

{\bf Remark.} The fact that only involutions may have non-zero
characters is known (see Corollary 1.10 in \cite{FJ2}, where
a combinatorial explanation of this fact is given). Our proof is
close to a similar proof from \cite{CT} where the case $G = S_n$
was considered.

Now we are ready to compute the characters for the classical
series.

Let us start with $A_{n-1}$, corresponding to the symmetric
group $G = S_n$ acting by permutations on the configuration space
$\mathcal M_n$ of $n$ distinct points of the complex plane. Any
involution $\sigma$ in that case is conjugated to one of the
involutions of the form $\pi_k = (12)(34)\dots(2k\!-\!1,2k), k
=1,\dots,[n/2],$ where $(ij)$ denote the  transposition of 
%%% I omitted ``elementary''
$i$ and $j.$ The fixed set $F_k$ for the action $\bar \pi_k$
consists of the points $(z_1,\bar z_1, z_2,\bar z_2,\dots, z_k,
\bar z_k, x_1,\dots,x_{n-2k}),$ where $z_i \neq \bar z_i$ (i.e. $z_i$
are not real), $x_i\in\mathbb R$,
$z_i \neq z_j, z_i \neq \bar z_j$ if $i \neq j$ and
$x_p \neq x_q$ if $p \neq q.$

It is easy to see that topologically all connected components are
the same and equivalent to the product $\mathcal M_k \times {\mathbb
R}^{n-2k},$ where $\mathcal M_k$ is the standard configuration space
of $k$ different points in the complex plane. Because the Euler
characteristic of $\mathcal M_k$ is zero for $k >1$ we conclude that
only the simple reflection $s = \pi_1 = (12)$ has non-zero character.
In that case we have $2(n-2)!$ contractible connected components,
so $\chi(s) = 2(n-2)!.$ We can formulate this fact in the
following form.

\begin{proposition} (G. Lehrer \cite{L1}) \label{2}
As $S_n$ module with respect to the standard action of $S_n$ on
$\mathcal M_n$ $$H^*(\mathcal M_n) = 2 Ind_{<s>}^{S_n}(1),$$ where
$Ind_{<s>}^{S_n}(1)$ denote the induced representation of $S_n$
from the trivial representation of the subgroup ${\mathbb Z}_2$
generated by the simple reflection $s = (12).$
\end{proposition}

Indeed the characters for the induced representations $Ind_H^G(1)$
from the trivial representation of the subgroup $H \subset G$ are
given by the formula
\begin{equation}
\label{ind} \chi(g) = \frac{1}{|H|} \sum_{h \sim g, h \in H}
|C(h)|,
\end{equation}
where $h \sim g$ means that $h$ is conjugate to $g$ in $G$, $C(h)$
is the centraliser of $h$: $C(h) = \{g \in G: gh = hg \}$. The
centraliser of the simple reflection $s$ is $C(s) = {\mathbb Z}_2
\times S_{n-2},$ where ${\mathbb Z}_2$ is generated by $s = (12)$ and
$S_{n-2}$ is the subgroup of $S_n$ permuting the numbers
$3,4,\dots,n.$ According to the formula (\ref{ind}) only identity
$e$ and $g \sim s$ have non-zero characters in  $2
Ind_{<s>}^{S_n}(1).$ Moreover $\chi(e) = |S_n| = n!$ and $\chi(s)
= |C(s)| = 2 (n-2)!$. Comparing this with our previous analysis we
have the proposition.

Now let us consider the twisted action of the symmetric group
$S_n$ on $\mathcal M_n.$ From the previous analysis of the fixed sets
$F_{\bar g}$ and Proposition 1 we deduce the following

\begin{proposition} \label{3}
The twisted action of $S_n$ on  $\mathcal M_n$ induces the regular
representation on the total cohomology $H^*(\mathcal M_n)$.
\end{proposition}

Consider now the case $B_n$. In the standard geometric
realisation any involution is conjugate either to the diagonal
form $\tau_k = diag (-1,\dots,-1,1,\dots,1) = (-I_k) \oplus I_{n-k}$ or
to the form when we have $m$ additional elementary transposition
blocks $P$: $\sigma_{m,k} = P \oplus P\dots\oplus P \oplus (-I_k)
\oplus I_{n-k-2m}.$ A simple analysis similar to the previous case
shows that in the last case the fixed set consists of several
connected components each topologically equivalent to $\mathcal M_m
\times {\mathbb R}^{n-2m}.$ So if the number of transpositions $m$ is
bigger than $1$ the Euler characteristic of this set and thus the
corresponding character $\chi(\sigma_{m,k})$ is zero. One can
check that in the remaining cases $\sigma_k = \sigma_{1,k}$ and
$\tau_k$ all the components are contractible and the action of the
corresponding centralisers $C(\sigma)$ on them is effective. Thus
if we define the set $X_{B_n}$ as the union of the $2n$
involutions $\sigma_k, k = 0,\dots,n-2$ and $\tau_l, l = 0,1,\dots,n$ we
have the following

\begin{proposition} \label{4}
\begin{equation}
\label{Bn}
H^*(\mathcal M_{B_n}) = \sum_{\sigma \in X_{B_n}} (2
Ind_{<\sigma>}^{B_n}(1) - \rho)
\end{equation} where $\rho$ is the
regular representation of $B_n$. In the twisted case one should
take the sum over the subset of $X_{B_n}$ consisting of the
involutions which are even elements of $G$.
\end{proposition}

For the standard action this result was first obtained by 
Lehrer in \cite{L2}.

In the $D_n$ case the results are slightly different for odd and
even $n$.

For an odd $n=2m+1$ any involution is conjugate to one of the following
involutions: $\sigma_{2k}, k =0,1,\dots,m-1$
and $\tau_{2l}, l = 0, 1,\dots, m $ (with the same notations as in the
previous case). 

The fixed set of $\sigma_{2k}$ consists of the points $(z,\bar z,
ix_1,\dots,ix_{2k}, y_1,\dots,y_{n-2k-2}),$ where $z \neq \pm \bar z$
(i.e. $z$ is not real or pure imaginary), $x_i,y_i\in\mathbb R$,
$x_p \neq \pm x_q$, $y_p
\neq \pm y_q$ if $p \neq q.$ The number of connected
components of the fixed set in that case is equal to the order of the group $D_2
\times D_{2k} \times D_{n-2k-2}$ which acts 
on this set. However if $k >0$ we have to exclude from this set also the subspaces
of codimension 2 given by the equation $x_j = y_l = 0$ for all $j=1,\dots, 2k$ and $l= 1,\dots, n-2k-2.$
One can check that as a result the Euler characteristic of all connected components
(and thus of the whole fixed set) in that case is $0.$ Similar arguments show that 
the same is true for the fixed set of the involutions $\tau_{2l}.$

In the only remaining case of simple reflection  $\sigma_0$ all the connected components of the corresponding fixed set are contractible. The centraliser $C(\sigma_0)$
coincides with the subgroup $D_2 \times D_{n-2},$ which acts freely and transitively
on these components.

For the groups of type $D_n$ with even $n=2m$ the involutions $\sigma_0, \sigma_{2m-2}$
and central element $\tau_{2m} = -I_n$ have the centralisers whose
order is equal to the number of connected components of the corresponding fixed sets,
which are all contractible in this case.
The Euler characteristic of the fixed sets of all other involutions can be shown to be zero
in the same way as in the odd case.

Thus we have the following
\begin{proposition} \label{5}
For the Coxeter group $G$ of type $D_{n}$ with the standard action we have
\begin{equation}
\label{Dodd} H^*(\mathcal M_G) = 2 Ind_{<\sigma_0>}^G(1) \end{equation}
in the odd case $n = 2m+1$ and 
\begin{equation}
\label{Deven} H^*(\mathcal M_G) = \rho + (2 Ind_{<\sigma_0>}^G(1) -
\rho) + (2 Ind_{<\sigma_{2m-2}>}^G(1)-\rho) + (2
Ind_{<\tau_{2m}>}^G(1) - \rho) 
\end{equation}
in the even case $n=2m.$
For the twisted action the formulas are respectively
\begin{equation}
\label{tDodd} H_{\epsilon}^*(\mathcal M_G) = \rho 
\end{equation}
in the odd case and
\begin{equation}
\label{tDeven} H_{\epsilon}^*(\mathcal M_G) = 2 Ind_{<\tau_{2m}>}^G(1)  
\end{equation}
in the even case.
\end{proposition}

In order to investigate the general case and to understand the
nature of the exceptional set $X_G$ we will need the general
results about involutions in the Coxeter groups summarised in the
next section.

\section*{Involutions in the Coxeter groups and their
centralisers.}

We start with the description of the conjugacy classes of
involutions in general Coxeter group due to Richardson \cite{R}.
We refer to Bourbaki \cite{Bourbaki} or Humphreys \cite{Hum} for
the standard definitions.

Let $G$ be a Coxeter group, $\Gamma$ be its Coxeter graph. Recall
that the vertices of $\Gamma$ correspond to the generating
reflections, which form the set $S.$ We will identify $G$ with its
geometric realisation in the Euclidean vector space $V,$ in which
$s \in S$ acts as the reflection with respect to a hyperplane
orthogonal to the corresponding root $e_s.$

Let $J$ be a subset of $S,$ $G_J$ be the subgroup $G$ generated by
$J$ (such a subgroup is called {\it parabolic}). Let also $J^* =
\{e_s, s \in J \}$ and $V_J$ be the subspace of $V$ generated by
$J^*$. Following Richardson we say that $J$ satisfies the {\it
$(-1)$-condition} if $G_J$ contains an element $\sigma_J$ which
acts on $V_J$ as $-Id.$ This element, being in $G_J$, acts
as the identity on the orthogonal complement of $V_J$. %%%added explanation
Thus $\sigma_J$ is uniquely determined and is an involution.

\begin{proposition} (Richardson \cite{R})
 Let $\sigma$ be any involution in $G$. Then
\begin{enumerate}
\item[(i)]
there exists a subset $J \subset S$ satisfying the $(-1)$-condition
such that $\sigma$ is conjugate to $\sigma_J$;
\item[(ii)]
%%% the involutions $\sigma_J$ and $\sigma_K$ are conjugate in $G$ if
the involutions $\sigma_J$ and $\sigma_K$ are conjugate in $G$ if
and only if $g(J^*) = K^*$ for some $g\in G$.
\end{enumerate}
\end{proposition}

Richardson gave also an algorithm for testing $G$-equivalence
based on the results by Howlett \cite{H} and Deodhar \cite{D}.

For an irreducible Coxeter group $G$ the whole set $J=S$ satisfies
the $(-1)$-condition only in the following cases: $$A_1, B_n, D_{2n},
E_7, E_8, F_4, H_3, H_4, I_2(2n).$$ This means that each connected
component of the Coxeter graph corresponding to the group $G_J$
for any $J$ satisfying the $(-1)$-condition must be of that form. In
particular we see that the components of type $A_n$ with $n>1$ are
forbidden, which imposes strong restrictions on such subsets $J.$
%New
In the Appendix we give a subset $J$ for each conjugacy class of special
involutions.

To describe the centralisers of the involutions $\sigma_J$ we can
use the results by Howlett \cite{H} who described the normalisers
of the parabolic subgroups in the general Coxeter group. Indeed we
have the following

\begin{proposition}
The centraliser $C(\sigma_J)$ coincides with the normaliser
$N(G_J)$ of the corresponding parabolic subgroup $G_J.$
\end{proposition}

Recall that the {\it normaliser} $N(H)$ of a subgroup $H \subset
G$ consists of the elements $g \in G$ such that $g H g^{-1} = H.$
Consider the orthogonal splitting $V = V_J \oplus V_J^{\perp},$
corresponding to the spectral decomposition of $\sigma_J.$
Obviously the normaliser $N(G_J)$ preserves this splitting and
thus $N(G_J)$ is a subgroup of the centraliser $C(\sigma_J).$ To
show the opposite inclusion take any $g \in C(\sigma_J)$ then
$g(V_J) = V_J$ and therefore $g(e_s)$ belongs to the intersection
of the root system $R$ with the subspace $V_J.$ But according to
\cite{R} (Proposition 1.10) this intersection coincides with the
root system of the group $G_J.$ This means that $g s g^{-1}$
belong to $G_J$ for any generating reflection $s \in J$ and thus
$g$ belongs to the normaliser $N(G_J).$

So the question now is only what are the special involutions which
form the set $X_G.$

\section*{Special involutions in the Coxeter groups.}

Let $\sigma$ be any involution in the Coxeter group $G,$ $V = V_1
\oplus V_2$ is the corresponding spectral decomposition of its
geometric realisation, where $V_1$ and $V_2$ are the eigenspaces
with the eigenvalues $-1$ and $1$, respectively. According to the
previous section one can assume that $V_1 = V_J$ for some subset
$J \in S$.

Consider the intersections $R_1$ and $R_2$ of the root system $R$
of $G$ with the subspaces $V_1$ and $V_2$. According to \cite{R},
\cite{Steinberg} they are the root systems of the Coxeter
subgroups $G_1$ and $G_2,$ where $G_1$ is the corresponding
parabolic subgroup $G_J$ and $G_2$ is the subgroup of $G$
consisting of the elements fixing the subspace $V_1.$

{\bf Definition.} We will call the involution $\sigma$ {\it
special} if for any root from $R$ at least one of its projections
onto $V_1$ and $V_2$ is proportional to a root from $R_1$ or
$R_2.$

In particular, the identity and simple reflections are special involutions for all
Coxeter groups. 

The following Proposition explains the importance of this notion
for our problem.

\begin{proposition}
For any special involution $\sigma$ all the connected components
of the fixed set $F_{\bar \sigma}$ are contractible. Their number
$N_{\sigma} = |G_1| |G_2|$ is equal to the product of the orders
of the corresponding groups $G_1$ and $G_2.$
\end{proposition}

{\bf Proof.} The fixed set $F_{\bar \sigma}$ is the subspace $V_2
\oplus i V_1$ minus the intersections with the reflection
hyperplanes. The hyperplanes corresponding to the roots from $R_1$
and $R_2$ split this subspace into $|G_1| |G_2|$ contractible
connected components. The meaning of the condition on the special
involutions is to make all other intersections redundant. Indeed
the intersection with the hyperplane $(\alpha, z) = 0$ is
equivalent to two conditions: $(\alpha_1,x) = 0$ and $(\alpha_2,
y) =0$ where $\alpha_1$ and $\alpha_2$ are the projections of the
root $\alpha$ on $V_1$ and $V_2$ respectively, $x \in V_1, y \in
V_2.$ If at least one of the vectors $\alpha_i$ is proportional to
some root from $R_i$ this intersection will be contained in the
hyperplanes already considered.

\begin{corollary} %%% reformulated
The character of the action of any special involution $\sigma$
on $H^*(\mathcal M_G)$ is non-zero and equal to $|G_1| |G_2|.$
\end{corollary}

It turns out that the converse statement is also true.

\begin{proposition}
If the character $\chi(g) \neq 0$ in $H^*(\mathcal M_G)$ then $g$ is a
special involution.
\end{proposition}

Unfortunately the only proof we have is case by case check. For
the classical series this follows from the fact that the special
involutions in that case are exactly those described in the first
section of this paper. For the exceptional Weyl groups one can
also use the results of \cite{FJ2}.

The list of the special involutions (up to a conjugation) for all irreducible Coxeter groups is given in the Introduction (see also the Appendix for the corresponding Richardson's graphs). It is easy to check that they all have the following property which is very important for us.

\begin{proposition}
The centraliser $C(\sigma)$ of any special involution $\sigma$ coincides with the product
of the corresponding Coxeter subgroups $G_1 \times G_2.$ 
\end{proposition}

Obviously $G_1 \times G_2$ is a subgroup of $C(\sigma)$ but the fact that $C(\sigma) = G_1 \times G_2$ is not true for a general involution (an example is any of the involutions $\sigma_{2k}$
with $k >0$ in $D_{2m+1}$-case, for which $G_1 \times G_2$ is a subgroup of $C(\sigma)$ of index 2).

Summarising we have the following main

\begin{theorem}
Let $X_G$ be the set of conjugacy classes of special involutions
in the Coxeter group $G$. Then the total cohomology
$H^*(\mathcal M_G)$ as $G$-module with respect to the standard action
of $G$ on $\mathcal M_G$ can be represented in the form
\begin{equation}
H^*(\mathcal M_G) = \sum_{\sigma \in X_G} (2
Ind_{<\sigma>}^G(1) - \rho).
\end{equation}
For the twisted action one should replace in this formula $X_G$ by
its subset consisting of the conjugacy classes of even special involutions.
\end{theorem}

In particular, we see that $H_{\epsilon}^*(\mathcal M_G)$ is the regular representation only for the symmetric group $S_n$, the Weyl groups of type $D_{2m+1}$, $E_6$
and dihedral groups $I_2(2k+1)$. 

Applying Frobenius reciprocity formula for the characters of the induced representations
we have the following corollary. Let $W$ be an irreducible representation of $G,$
$W^\pm(\sigma)=\{v\in W\colon \sigma v=\pm v\}$ 
be the eigenspaces of the involution $\sigma$ with eigenvalues $\pm 1$ respectively.

\begin{corollary}
The multiplicity $m(W)$ in the decomposition of $H^*(\mathcal M_G)$ is equal to
$$
m(W) = \sum_{\sigma\in X_G}(\mathrm{dim}\,W^+(\sigma)-\mathrm{dim}\,W^-(\sigma)).
$$
In particular, for the trivial and alternating representations
we have respectively
$$
m(1) = |X_G|
$$
and 
$$
m(\epsilon) = |X^{\epsilon}_G| - |X^{o}_G|,
$$
where $X^{o}_G$ is the subset of $X_G$ consisting of odd special involutions and 
$|X|$ denotes the number of elements in the set $X$. 
\end{corollary}

In other words the number of special involutions (up to a conjugacy) is equal to the 
number of invariant cohomology classes in $H^*(\mathcal M_G).$ 
The corresponding numbers are given in the table below.

\bigskip

\begin{tabular} {|c|c|c|c|c|c|c|c|c|c|c|c|}
\hline
\multicolumn{12}{|c|}{\textbf{Numbers $|X_G|$ of conjugacy classes of special involutions}}\\
\hline
$A_n$ & $B_n$ & $D_n$, $n$ odd & $D_n$, $n$ even & $E_6$ & $E_7$& $E_8$ & $F_4$ & $H_3$ &$H_4$ & $I_2(n)$, $n$ odd &$I_2(n)$, $n$ even \\
\hline
2&$2n$&2&4&2&4&4&8&4&4&2&4\\
\hline
\end{tabular}

\bigskip

A simple analysis of this list leads to the following

\begin{proposition}
%\label{anti}
For any irreducible Coxeter group $G$ the numbers of even and odd involutions in the special set $X_G$ are equal:  
$$|X^{o}_G| = |X^{\epsilon}_G| = \frac{1}{2} |X_G|. $$
The multiplicity $m(\epsilon)$ of the alternating representation in $H^*(\mathcal M_G)$ is 
zero for all $G$.
\end{proposition}

In other words the alternating representation never appears in the cohomology of $\mathcal M_G.$
For the symmetric group the fact that  $m(\epsilon) = 0$ was conjectured by Stanley \cite{Stanley}
(in a different, combinatorial terms) and proved by  Hanlon in \cite{Hanlon}. The general case was settled by Lehrer in \cite{L3}.  We give some explanations of this remarkable fact in the next section.

\section*{Special involutions and cohomology of Brieskorn's braid groups.}

Consider now the corresponding quotient
space $\Sigma_G = \mathcal M_G / G.$
The fundamental group of this space is called {\it Brieskorn's braid group} and denoted as $B_G.$
For the space $\mathcal M_G$ it is called {\it Brieskorn's pure braid group} and denoted as $P_G.$

It is known after Brieskorn and Deligne that $\mathcal M_G$ and  $\Sigma_G$
are Eilenberg--Mac Lane spaces $K(\pi, 1)$ for the pure and usual braid groups
$P_G$ and $B_G$ respectively (see \cite{Brieskorn}, \cite{Deligne}), so the cohomology of these spaces
can be interpreted as cohomology of the corresponding braid groups.
The investigation of these cohomology was initiated by Arnol'd in \cite{Arnold}, \cite{Arnold2}
who considered the case of the symmetric group. The rational (or complex) coefficients cohomology of the braid groups related to arbitrary Coxeter group were computed by Brieskorn \cite{Brieskorn}.
Later Orlik and Solomon \cite{OS} found an elegant description of the de Rham
cohomology as the algebra generated by the differential forms $\omega_{\alpha} = d \, \log \alpha$, $\alpha\in R$, with explicitly given relations (defining the so-called {\it Orlik--Solomon algebra} $\mathcal{A}_G$).

Since $H^*(\Sigma_G, \mathbb{C})$ is simply the $G$-invariant part of $H^*(\mathcal M_G, \mathbb{C})$, it follows 
from the previous section that the total Betti number of $\Sigma_G$, which is the dimension of 
$H^*(\Sigma_G, \mathbb{C})$, is equal to $|X_G|.$ The relation between special involutions in $G$ and the cohomology $H^*(\Sigma_G, \mathbb{C})$ can be described more precisely in the following way.

\begin{proposition}
%\label{ups}
The Poincar\'e polynomial $P (\Sigma_G, t)$ of the cohomology of the Brieskorn's braid group $B_G$
has the form
\begin{equation}
\label{pinv}
P (\Sigma_G, t) = \sum_{\sigma\in X_G} t^{\mathrm{dim}\,V^-(\sigma)},
\end{equation}
where $V$ is the geometric representation of $G$ and $V^-(\sigma)$ is the $(-1)$-eigenspace
of $\sigma$ in this representation.
\end{proposition}

Explicitly, these polynomials are

$A_n: \quad P = 1 + t$

$B_n: \quad P = 1 + 2t + 2t^2 + \dots + 2t^{n-1} + t^n = (1+ t ) (1+ t+ t^2 +\dots + t^{n-1})$

$D_n$, $n$ odd: $\quad P = 1 + t$

$D_n$, $n$ even: $\quad P = 1 + t + t^{n-1} + t^n = (1 + t) (1 + t^{n-1})$

$E_6: \quad P = 1 +t$

$E_7: \quad P = 1 +t + t^6 + t^7 = (1 + t) (1 + t ^6)$

$E_8: \quad P = 1 +t + t^7 + t^8 = (1 + t) (1 + t ^7)$

$F_4: \quad P = 1 +2t + 2t^2 + 2t^3 + t^4= (1 + t) (1 + t + t^2 + t^3)$

$H_3: \quad P = 1 +t +t^2 + t^3 = (1+ t ) (1 + t^2)$

$H_4: \quad P = 1 +t +t^3 + t^4 = (1+ t ) (1 + t^3)$

$I_2(n)$, $n$ odd: $\quad P = 1+ t$

$I_2(n)$, $n$ even: $\quad P = 1+ 2t + t^2 = (1+t)^2.$

The simplest proof of this proposition is by comparison of the Brieskorn formulas from \cite{Brieskorn}
and our list of simple involutions.  A more satisfactory proof may be found in the following way in terms of the corresponding Orlik--Solomon algebra $\mathcal{A}_G \approx H^*(\mathcal M_G)$.

Let $\sigma$ be a special involution and $S_{\sigma}$ be a set of the simple roots in the Coxeter subsystem $R_1,$  which is the set of roots $\alpha \in R$
such that $\sigma \alpha = - \alpha $ (see the previous section).

{\bf Conjecture.} {\it For any special involution $\sigma$ the symmetrisation of the product
$\Omega_{\sigma} = \prod_{\alpha \in S_{\sigma}} \omega_{ \alpha}$
by the action of the group $G$ is a non-zero element in $\mathcal{A}_G.$ Any $G$-invariant element in this algebra is a linear combination of such elements.}

We have checked this for all Coxeter groups except $E_7, E_8, F_4, H_3, H_4,$
but the proof is by straightforward calculation with the use of the Orlik--Solomon relations.
It would be nice to find a more conceptual proof for all Coxeter groups.

Notice that in $S_n$ case our claim reduces to Arnold's result saying that the symmetrisation of any element from the cohomology of the pure braid group of degree more than 1 is equal to zero and the only $G$-invariant element of degree 1 is given by the differential form $\sum_{i \neq j} d \, \log (z_i - z_j)$
(see Corollary 6 in \cite{Arnold}).

Notice that as a corollary of Proposition 11 we have that the anti-symmetrisation of any element
in Orlik--Solomon algebra is zero. This is also related to the fact that all the polynomials 
$P (\Sigma_G, t)$ are divisible by $t+1.$ Indeed from the previous section we know that 
$m(\epsilon) = |X^{\epsilon}_G| - |X^{o}_G|$ which due to Proposition 12 is equal to $P (\Sigma_G, -1)$
and thus is zero. 

Although the fact that $P (\Sigma_G, t)$ is divisible by $t+1$ is transparent from the list of all these polynomials given above we will give here an independent topological explanation.

Namely consider the {\it projectivisation} $P\mathcal M_G$ of the space $\mathcal M_G.$
We have a diffeomorphism
$$\mathcal M_G \approx  P\mathcal M_G \times {\mathbb{C}}^*.$$
Indeed for any root $\alpha$ the map $x \rightarrow ([x], (\alpha, x))$
establish such a diffeomorphism. Notice that it is compatible with the action of the
group $G$ if we assume that the action on ${\mathbb{C}}^*$ is trivial.

On the cohomology level we have an isomorphism
$$H^*(\mathcal M_G) =  H^*(P\mathcal M_G) \times H^*({\mathbb{C}}^*),$$
which immediately implies the following result.

Let $W$ be any irreducible representation of $G$ and $P_W (t)$ be the Poincar\'e polynomials of
the multiplicities of $W$ in $H^k(\mathcal M_G).$ An explicit description of these polynomials is one of the most interesting open problems in this area.
For the trivial representation the corresponding polynomial coincides with the Poincar\'e polynomial $P (\Sigma_G, t)$ according to Proposition 12.

\begin{proposition}
All irreducible representations of the Coxeter group $G$ appear in $H^*(\mathcal M_G)$ in pairs in degrees 
differing by 1. The corresponding polynomials $P_W(t)$ are divisible by $t+1.$
\end{proposition}

In particular, this implies that $P (\Sigma_G, -1) = 0$ and (modulo Propositions 11 and 12)
the fact that the alternating representation never appears in the cohomology of generalised pure braid
groups.

\section*{Concluding remarks.}

We have shown that the action of the Coxeter groups $G$ on the total cohomology of the corresponding
complement space $\mathcal M_G$ can be described in a very simple way (\ref{F1}) in terms of the special involutions. Although one might expect a formula like that knowing that only involutions may have non-zero characters in this representation, there are two facts which we found surprising:

1) the involutions with non-zero characters admit a very simple geometric characterisation;

2) the character of the action of such an involution $\sigma$ on $H^*(\mathcal M_G)$ is exactly the order of the corresponding centraliser $C(\sigma).$ 

% 3) the numbers of even and odd involutions in the special set $X_G$ are always equal.
%New
Special involutions also appear to play a key role in the description of the multiplicative structure of the cohomology of the Brieskorn braid group: if the conjecture formulated in the last section is true, a basis of representatives of the cohomology in the Orlik--Solomon algebra can be given in terms of special invlutions.

The nature of the set $X_G$ also needs better understanding. In particular, it may be worth looking at the geometry of the subgraphs $\Gamma_J$
in the Coxeter graph $\Gamma$ corresponding to the special
involutions (see the Appendix for a complete list of these subgraphs). We would like to note that in all the cases except $D_4$ $\Gamma_J$ consist of at most two
connected components which together with the $(-1)$-condition implies
very strong restrictions on $J$. The only exception is the graph of the involution $-s$ in $D_4,$ which consists of three components.

It would be interesting to generalise our results for the space
$\mathcal M_G({\mathbb R}^N)$ which is $V \otimes {\mathbb R}^N$ without
corresponding root subspaces $\Pi_{\alpha} \otimes {\mathbb R}^N,
\alpha \in R$ (our case corresponds to $N=2$). For the symmetric group $G = S_n$ 
this is the configuration space of $n$ distinct points in ${\mathbb R}^N,$ this case was investigated by Cohen and Taylor in \cite{CT}. In this
relation we would like to mention the recent very interesting papers
\cite{Atiyah, A2, AB, FY}. In particular, the idea of Atiyah
\cite{A2} to use equivariant cohomology may be the clue to
understanding the graded structure of the $G$-module
$H^*(\mathcal M_G).$

Another obvious generalisation to look at is the case of complex reflection groups.

\section*{ Acknowledgements}

We are grateful to  Corrado De Concini, Eva Feichtner, Gus Lehrer, Guido Mislin, Richard Stanley, Victor Vassiliev and especially to Misha Farber, who attracted our attention to the
generalised Lefschetz formula and to the Cohen--Taylor paper \cite{CT},
for useful discussions and comments.

One of us (APV) is grateful to the Forschungsinstitut f\"ur
Mathematik (ETH, Zurich) for the hospitality in the spring 2003
when this work had started.

\newpage
%Coxeter graphs
\setlength\unitlength{0.8cm}

%An
\newcommand{\An}{
\begin{picture}(8,0.8)(0,1.8)
\put(1,2){\circle*{0.2}}
\put(1.1,2){\line(1,0){.8}}
\put(2,2){\circle{0.2}}
\put(2.1,2){\line(1,0){.8}}
\put(3,2){\circle{0.2}}
\put(3.1,2){\line(1,0){.2}}
\put(3.55,2){\circle*{.02}}
\put(3.75,2){\circle*{.02}}
\put(3.95,2){\circle*{.02}}
\put(4.2,2){\line(1,0){.2}}
\put(4.5,2){\circle{0.2}}
\put(4.6,2){\line(1,0){.8}}
\put(5.5,2){\circle{0.2}}
\end{picture}
}

% Bn
\newcommand{\Bna}
{\begin{picture}(8,0.8)(0.5,1.8)
\put(1,2){\circle*{0.2}}
\put(1.1,2){\line(1,0){.6}}
\put(1.8,2){\circle{0.2}}
\put(1.9,2){\line(1,0){.6}}
\put(2.6,2){\circle{0.2}}
\put(2.7,2){\line(1,0){.2}}
\put(3.15,2){\circle*{.02}}
\put(3.35,2){\circle*{.02}}
\put(3.55,2){\circle*{.02}}
\put(3.8,2){\line(1,0){.2}}
\put(4.1,2){\circle{0.2}}
\put(4.2,2){\line(1,0){.6}}
\put(4.9,2){\circle*{0.2}}
\put(5,2){\line(1,0){.6}}
\put(5.7,2){\circle*{0.2}}
\put(5.8,2){\line(1,0){.2}}
\put(6.25,2){\circle*{.02}}
\put(6.45,2){\circle*{.02}}
\put(6.65,2){\circle*{.02}}
\put(6.9,2){\line(1,0){.2}}
\put(7.2,2){\circle*{0.2}}
\put(7.25,2.05){\line(1,0){.7}}
\put(7.25,1.95){\line(1,0){.7}}
\put(8,2){\circle*{0.2}}
\end{picture}
}
\newcommand{\Bnb}
{\begin{picture}(8,0.8)(0.5,1.8)
\put(1,2){\circle{0.2}}
\put(1.1,2){\line(1,0){.6}}
\put(1.8,2){\circle{0.2}}
\put(1.9,2){\line(1,0){.6}}
\put(2.6,2){\circle{0.2}}
\put(2.7,2){\line(1,0){.2}}
\put(3.15,2){\circle*{.02}}
\put(3.35,2){\circle*{.02}}
\put(3.55,2){\circle*{.02}}
\put(3.8,2){\line(1,0){.2}}
\put(4.1,2){\circle{0.2}}
\put(4.2,2){\line(1,0){.6}}
\put(4.9,2){\circle*{0.2}}
\put(5,2){\line(1,0){.6}}
\put(5.7,2){\circle*{0.2}}
\put(5.8,2){\line(1,0){.2}}
\put(6.25,2){\circle*{.02}}
\put(6.45,2){\circle*{.02}}
\put(6.65,2){\circle*{.02}}
\put(6.9,2){\line(1,0){.2}}
\put(7.2,2){\circle*{0.2}}
\put(7.25,2.05){\line(1,0){.7}}
\put(7.25,1.95){\line(1,0){.7}}
\put(8,2){\circle*{0.2}}
\end{picture}
}
%H3
\newcommand{\Hthreea}{
\begin{picture}(8,0.8)(0,1.8)
\put(1,2){\circle*{0.2}}
\put(1.4,2.1){5}
\put(1.1,2){\line(1,0){.8}}
\put(2,2){\circle{0.2}}
\put(2.1,2){\line(1,0){.8}}
\put(3,2){\circle{0.2}}
\end{picture}
}
\newcommand{\Hthreeb}{
\begin{picture}(8,0.8)(0,1.8)
\put(1,2){\circle*{0.2}}
\put(1.4,2.1){5}
\put(1.1,2){\line(1,0){.8}}
\put(2,2){\circle{0.2}}
\put(2.1,2){\line(1,0){.8}}
\put(3,2){\circle*{0.2}}
\end{picture}
}
\newcommand{\Hthreec}{
\begin{picture}(8,0.8)(0,1.8)
\put(1,2){\circle*{0.2}}
\put(1.4,2.1){5}
\put(1.1,2){\line(1,0){.8}}
\put(2,2){\circle*{0.2}}
\put(2.1,2){\line(1,0){.8}}
\put(3,2){\circle*{0.2}}
\end{picture}
}
%I2(m)
\newcommand{\Itwoa}{
\begin{picture}(8,0.8)(0,1.8)
\put(1,2){\circle*{0.2}}
\put(1.3,2.1){$m$}
\put(1.1,2){\line(1,0){.8}}
\put(2,2){\circle{0.2}}
\end{picture}
}
\newcommand{\Itwob}{           
\begin{picture}(8,0.8)(0,1.8)
\put(1,2){\circle{0.2}}
\put(1.3,2.1){$m$}
\put(1.1,2){\line(1,0){.8}}
\put(2,2){\circle*{0.2}}
\end{picture}
}
\newcommand{\Itwoc}           
{
\begin{picture}(8,0.8)(0,1.8)
\put(1.3,2.1){$m$}
\put(1,2){\circle*{0.2}}
\put(1.1,2){\line(1,0){.8}}
\put(2,2){\circle*{0.2}}
\end{picture}
}
%E6
\newcommand{\Esix}{
\begin{picture}(8,1.5)(0,1.8)
\put(1,2){\circle*{0.2}}
\put(1.1,2){\line(1,0){.8}}
\put(2,2){\circle{0.2}}
\put(2.1,2){\line(1,0){.8}}
\put(3,2){\circle{0.2}}
\put(3.1,2){\line(1,0){.8}}
\put(4,2){\circle{0.2}}
\put(4.1,2){\line(1,0){.8}}
\put(5  ,2){\circle{0.2}}
\put(3  ,2.1){\line(0,1){.8}}
\put(3  ,3){\circle{0.2}}
\end{picture}
}
%E7
\newcommand{\Esevena}{
\begin{picture}(8,1.5)(0,1.8)
\put(1,2){\circle*{0.2}}
\put(1.1,2){\line(1,0){.8}}
\put(2,2){\circle{0.2}}
\put(2.1,2){\line(1,0){.8}}
\put(3,2){\circle{0.2}}
\put(3.1,2){\line(1,0){.8}}
\put(4,2){\circle{0.2}}
\put(4.1,2){\line(1,0){.8}}
\put(5  ,2){\circle{0.2}}
\put(3  ,2.1){\line(0,1){.8}}
\put(3  ,3){\circle{0.2}}
\put(5.1,2){\line(1,0){.8}}
\put(6,2){\circle{0.2}}
\end{picture}
}
\newcommand{\Esevenb}{
\begin{picture}(8,1.5)(0,1.8)
\put(1,2){\circle{0.2}}
\put(1.1,2){\line(1,0){.8}}
\put(2,2){\circle*{0.2}}
\put(2.1,2){\line(1,0){.8}}
\put(3,2){\circle*{0.2}}
\put(3.1,2){\line(1,0){.8}}
\put(4,2){\circle*{0.2}}
\put(4.1,2){\line(1,0){.8}}
\put(5  ,2){\circle*{0.2}}
\put(3  ,2.1){\line(0,1){.8}}
\put(3  ,3){\circle*{0.2}}
\put(5.1,2){\line(1,0){.8}}
\put(6,2){\circle*{0.2}}
\end{picture}
}
\newcommand{\Esevenc}{
\begin{picture}(8,1.5)(0,1.8)
\put(1,2){\circle*{0.2}}
\put(1.1,2){\line(1,0){.8}}
\put(2,2){\circle*{0.2}}
\put(2.1,2){\line(1,0){.8}}
\put(3,2){\circle*{0.2}}
\put(3.1,2){\line(1,0){.8}}
\put(4,2){\circle*{0.2}}
\put(4.1,2){\line(1,0){.8}}
\put(5  ,2){\circle*{0.2}}
\put(3  ,2.1){\line(0,1){.8}}
\put(3  ,3){\circle*{0.2}}
\put(5.1,2){\line(1,0){.8}}
\put(6,2){\circle*{0.2}}
\end{picture}
}
%E8
\newcommand{\Eeighta}{
\begin{picture}(8,1.5)(0,1.8)
\put(1,2){\circle*{0.2}}
\put(1.1,2){\line(1,0){.8}}
\put(2,2){\circle{0.2}}
\put(2.1,2){\line(1,0){.8}}
\put(3,2){\circle{0.2}}
\put(3.1,2){\line(1,0){.8}}
\put(4,2){\circle{0.2}}
\put(4.1,2){\line(1,0){.8}}
\put(5  ,2){\circle{0.2}}
\put(3  ,2.1){\line(0,1){.8}}
\put(3  ,3){\circle{0.2}}
\put(5.1,2){\line(1,0){.8}}
\put(6,2){\circle{0.2}}
\put(6.1,2){\line(1,0){.8}}
\put(7,2){\circle{0.2}}
\end{picture}
}

\newcommand{\Eeightb}{
\begin{picture}(8,1.5)(0,1.8)
\put(1,2){\circle*{0.2}}
\put(1.1,2){\line(1,0){.8}}
\put(2,2){\circle*{0.2}}
\put(2.1,2){\line(1,0){.8}}
\put(3,2){\circle*{0.2}}
\put(3.1,2){\line(1,0){.8}}
\put(4,2){\circle*{0.2}}
\put(4.1,2){\line(1,0){.8}}
\put(5  ,2){\circle*{0.2}}
\put(3  ,2.1){\line(0,1){.8}}
\put(3  ,3){\circle*{0.2}}
\put(5.1,2){\line(1,0){.8}}
\put(6,2){\circle*{0.2}}
\put(6.1,2){\line(1,0){.8}}
\put(7,2){\circle{0.2}}
\end{picture}
}

\newcommand{\Eeightc}{
\begin{picture}(8,1.5)(0,1.8)
\put(1,2){\circle*{0.2}}
\put(1.1,2){\line(1,0){.8}}
\put(2,2){\circle*{0.2}}
\put(2.1,2){\line(1,0){.8}}
\put(3,2){\circle*{0.2}}
\put(3.1,2){\line(1,0){.8}}
\put(4,2){\circle*{0.2}}
\put(4.1,2){\line(1,0){.8}}
\put(5  ,2){\circle*{0.2}}
\put(3  ,2.1){\line(0,1){.8}}
\put(3  ,3){\circle*{0.2}}
\put(5.1,2){\line(1,0){.8}}
\put(6,2){\circle*{0.2}}
\put(6.1,2){\line(1,0){.8}}
\put(7,2){\circle*{0.2}}
\end{picture}
}
%F4
\newcommand{\Ffoura}{
\begin{picture}(8,0.8)(0,1.8)
\put(1,  2){\circle*{0.2}}
\put(1.1,2){\line(1,0){.8}}
\put(2,  2){\circle{0.2}}
\put(2.1,2.05){\line(1,0){.8}}
\put(2.1,1.95){\line(1,0){.8}}
\put(3,  2){\circle{0.2}}
\put(3.1,2){\line(1,0){.8}}
\put(4,  2){\circle{0.2}}
\end{picture}
}
\newcommand{\Ffourb}{
\begin{picture}(8,0.8)(0,1.8)
\put(1,  2){\circle{0.2}}
\put(1.1,2){\line(1,0){.8}}
\put(2,  2){\circle{0.2}}
\put(2.1,2.05){\line(1,0){.8}}
\put(2.1,1.95){\line(1,0){.8}}
\put(3,  2){\circle{0.2}}
\put(3.1,2){\line(1,0){.8}}
\put(4,  2){\circle*{0.2}}
\end{picture}
}
\newcommand{\Ffourc}{
\begin{picture}(8,0.8)(0,1.8)
\put(1,  2){\circle{0.2}}
\put(1.1,2){\line(1,0){.8}}
\put(2,  2){\circle*{0.2}}
\put(2.1,2.05){\line(1,0){.8}}
\put(2.1,1.95){\line(1,0){.8}}
\put(3,  2){\circle*{0.2}}
\put(3.1,2){\line(1,0){.8}}
\put(4,  2){\circle{0.2}}
\end{picture}
}
\newcommand{\Ffourd}{
\begin{picture}(8,0.8)(0,1.8)
\put(1,  2){\circle*{0.2}}
\put(1.1,2){\line(1,0){.8}}
\put(2,  2){\circle{0.2}}
\put(2.1,2.05){\line(1,0){.8}}
\put(2.1,1.95){\line(1,0){.8}}
\put(3,  2){\circle*{0.2}}
\put(3.1,2){\line(1,0){.8}}
\put(4,  2){\circle{0.2}}
\end{picture}
}
\newcommand{\Ffoure}{
\begin{picture}(8,0.8)(0,1.8)
\put(1,  2){\circle{0.2}}
\put(1.1,2){\line(1,0){.8}}
\put(2,  2){\circle*{0.2}}
\put(2.1,2.05){\line(1,0){.8}}
\put(2.1,1.95){\line(1,0){.8}}
\put(3,  2){\circle*{0.2}}
\put(3.1,2){\line(1,0){.8}}
\put(4,  2){\circle*{0.2}}
\end{picture}
}
\newcommand{\Ffourf}{
\begin{picture}(8,0.8)(0,1.8)
\put(1,  2){\circle*{0.2}}
\put(1.1,2){\line(1,0){.8}}
\put(2,  2){\circle*{0.2}}
\put(2.1,2.05){\line(1,0){.8}}
\put(2.1,1.95){\line(1,0){.8}}
\put(3,  2){\circle*{0.2}}
\put(3.1,2){\line(1,0){.8}}
\put(4,  2){\circle{0.2}}
\end{picture}
}
\newcommand{\Ffourg}{
\begin{picture}(8,0.8)(0,1.8)
\put(1,  2){\circle*{0.2}}
\put(1.1,2){\line(1,0){.8}}
\put(2,  2){\circle*{0.2}}
\put(2.1,2.05){\line(1,0){.8}}
\put(2.1,1.95){\line(1,0){.8}}
\put(3,  2){\circle*{0.2}}
\put(3.1,2){\line(1,0){.8}}
\put(4,  2){\circle*{0.2}}
\end{picture}
}

%H4
\newcommand{\Hfoura}{
\begin{picture}(8,0.8)(0,1.8)
\put(1,  2){\circle*{0.2}}
\put(1.1,2){\line(1,0){.8}}
\put(1.4,2.1){5}
\put(2,  2){\circle{0.2}}
\put(2.1,2){\line(1,0){.8}}
\put(3,  2){\circle{0.2}}
\put(3.1,2){\line(1,0){.8}}
\put(4,  2){\circle{0.2}}
\end{picture}
}

\newcommand{\Hfourb}{
\begin{picture}(8,0.8)(0,1.8)
\put(1,  2){\circle*{0.2}}
\put(1.1,2){\line(1,0){.8}}
\put(1.4,2.1){5}
\put(2,  2){\circle*{0.2}}
\put(2.1,2){\line(1,0){.8}}
\put(3,  2){\circle*{0.2}}
\put(3.1,2){\line(1,0){.8}}
\put(4,  2){\circle{0.2}}
\end{picture}
}

\newcommand{\Hfourc}{
\begin{picture}(8,0.8)(0,1.8)
\put(1,  2){\circle*{0.2}}
\put(1.1,2){\line(1,0){.8}}
\put(1.4,2.1){5}
\put(2,  2){\circle*{0.2}}
\put(2.1,2){\line(1,0){.8}}
\put(3,  2){\circle*{0.2}}
\put(3.1,2){\line(1,0){.8}}
\put(4,  2){\circle*{0.2}}
\end{picture}
}
%Dn                                
\newcommand{\Dna}{
\begin{picture}(8,2)(0,1.8)
\put(1,2){\circle*{0.2}}
\put(1.1,2){\line(1,0){.8}}
\put(2,2){\circle{0.2}}
\put(2.1,2){\line(1,0){.8}}
\put(3,2){\circle{0.2}}
\put(3.1,2){\line(1,0){.2}}
\put(3.55,2){\circle*{.02}}
\put(3.75,2){\circle*{.02}}
\put(3.95,2){\circle*{.02}}
\put(4.2,2){\line(1,0){.2}}
\put(4.5,2){\circle{0.2}}
\put(4.6,2){\line(1,0){.8}}
\put(5.5,2){\circle{0.2}}
\put(5.5707,2.0707){\line(1,1){.859}}
\put(6.5,3){\circle{0.2}}
\put(5.5707,1.9293){\line(1,-1){.859}}
\put(6.5,1){\circle{0.2}}
\end{picture}
}                               
\newcommand{\Dnb}{
\begin{picture}(8,2.5)(0,1.8)
\put(1,2){\circle*{0.2}}
\put(1.1,2){\line(1,0){.8}}
\put(2,2){\circle{0.2}}
\put(2.1,2){\line(1,0){.8}}
\put(3,2){\circle*{0.2}}
\put(3.1,2){\line(1,0){.2}}
\put(3.55,2){\circle*{.02}}
\put(3.75,2){\circle*{.02}}
\put(3.95,2){\circle*{.02}}
\put(4.2,2){\line(1,0){.2}}
\put(4.5,2){\circle*{0.2}}
\put(4.6,2){\line(1,0){.8}}
\put(5.5,2){\circle*{0.2}}
\put(5.5707,2.0707){\line(1,1){.859}}
\put(6.5,3){\circle*{0.2}}
\put(5.5707,1.9293){\line(1,-1){.859}}
\put(6.5,1){\circle*{0.2}}
\end{picture}
}                               
\newcommand{\Dnc}{
\begin{picture}(8,2.5)(0,1.8)
\put(1,2){\circle*{0.2}}
\put(1.1,2){\line(1,0){.8}}
\put(2,2){\circle*{0.2}}
\put(2.1,2){\line(1,0){.8}}
\put(3,2){\circle*{0.2}}
\put(3.1,2){\line(1,0){.2}}
\put(3.55,2){\circle*{.02}}
\put(3.75,2){\circle*{.02}}
\put(3.95,2){\circle*{.02}}
\put(4.2,2){\line(1,0){.2}}
\put(4.5,2){\circle*{0.2}}
\put(4.6,2){\line(1,0){.8}}
\put(5.5,2){\circle*{0.2}}
\put(5.5707,2.0707){\line(1,1){.859}}
\put(6.5,3){\circle*{0.2}}
\put(5.5707,1.9293){\line(1,-1){.859}}
\put(6.5,1){\circle*{0.2}}
\end{picture}
}

\section*{Appendix: Richardson graphs of special involutions}
The following table lists the equivalence classes of nontrivial special involutions for all Coxeter groups. 
For each class we give a representative 
$\sigma_J$ associated with a subset $J$ of the set of nodes of the
Coxeter graph by the Richardson correspondence
(see Prop.~6):
the nodes in $J$, whose simple roots span the 
$-1$-eigenspace of $\sigma$, are coloured in 
black. The notations in the third column follow the description in the Introduction. 
Additionally to the special involutions listed 
here there is the identity, which is special for all groups and corresponds to a Coxeter graph with white nodes only.
The conventions for Coxeter graphs are those of
\cite{Bourbaki} Ch.~IV, \S 1, n${}^o$  9 
except that we draw a double edge
instead of an edge with label 4.
\bigskip

\begin{center}
\begin{tabular}{|l|l|l|}
\hline
$G$   &  $(\Gamma,J)$& $\sigma_J$\hspace{2cm} \\
\hline
$A_n$ & \An & $s$ \\
\hline
$B_n$ & \Bna & $\sigma_k$ \\
      & \Bnb & $\tau_l$ \\
\hline
$D_n$, $n$ odd & \Dna & $s$ \\
 & & \\
 & & \\
\hline
                  & \Dna & $s$ \\
$D_n$, $n$ even \hspace{0.4cm}    & \Dnb & $-s$ \\
                  & \Dnc & $-\mathit{Id}$ \\
 & & \\ & & \\
\hline
$E_6$ & \Esix & $ s$\\
\hline
      & \Esevena & $ s$\\
$E_7$ & \Esevenb & $ -s$\\
      & \Esevenc & $ -\mathit{Id}$\\
\hline
\end{tabular}

\begin{tabular}{|l|l|l|}
\hline
      & \Eeighta & $ s$\\
$E_8$ & \Eeightb & $ -s$\\
      & \Eeightc & $ -\mathit{Id}$\\
\hline
      & \Ffoura  & $s_1$\\
      & \Ffourb  & $s_2$\\
      & \Ffourc  & $diag(-1,-1,1,1) $\\
$F_4$ & \Ffourd  & $P_{12}\oplus diag(-1,1)$\\
      & \Ffoure  & $-s_1$\\
      & \Ffourf  & $-s_2$\\
      & \Ffourg  & $-\mathit{Id}$\\
\hline
      & \Hthreea  & $s$\\
$H_3$ & \Hthreeb  & $-s$\\
      & \Hthreec  & $-\mathit{Id}$\\
\hline
      & \Hfoura  & $s$\\
$H_4$ & \Hfourb  & $-s$\\
      & \Hfourc  & $-\mathit{Id}$\\
\hline
$I_2(m)$, $m$ odd & \Itwoa & $s$\\
\hline
                   & \Itwoa & $s_1$\\ 
$I_2(m)$, $m$ even & \Itwob & $s_2$\\ 
                   & \Itwoc & $-\mathit{Id}$\\ 
\hline

\end{tabular}
\end{center}

\end{document}